\documentclass{article}

\usepackage[noblocks]{authblk}

\usepackage{amssymb}
\usepackage{amsmath}
\usepackage{bbm}

\usepackage{amsthm}

\usepackage[T1]{fontenc}
\usepackage{lmodern}

\usepackage[pdftex]{graphicx}
\usepackage[captionskip=0mm,caption=false]{subfig}
\usepackage{float}
\graphicspath{{figures/}}

\usepackage{caption}
\usepackage{microtype}

\usepackage{verbatim}

\usepackage{placeins}

\usepackage[hyperfootnotes=false]{hyperref}

\usepackage{color}

\usepackage{bm}
\DeclareBoldMathCommand\boldE{\mathcal{E}}

\DeclareMathOperator{\E}{\mathcal{E}}

\DeclareMathOperator{\F}{\mathcal{F}}
\DeclareMathOperator{\Id}{\operatorname{Id}}
\renewcommand{\div}{\operatorname{div}}
\newcommand{\norm}[1]{\left\| #1 \right\|}
\newcommand{\iprod}[1]{\left\langle #1 \right\rangle}
\newcommand{\dn}[1]{\frac{\partial #1}{\partial \nu}}
\newcommand{\dnflat}[1]{\partial #1 / \partial \nu}
\newcommand{\vectwo}[2]{\ensuremath{(\begin{smallmatrix}#1\\#2\end{smallmatrix})}}
\def \Rm {\mathbb R}

\RequirePackage{ifthen}
\makeatletter
\newcommand{\logmessage}[1]{\@latex@warning{#1}}
\makeatother
\IfFileExists{../Bibliographie/Users_Guide.txt}{
  }{
  \IfFileExists{../../Bibliographie/Users_Guide.txt}{
    }{
    \IfFileExists{../../../Bibliographie/Users_Guide.txt}{
      }{
        \IfFileExists{../../../../Bibliographie/Users_Guide.txt}{
          }{
            \IfFileExists{../../../../../Bibliographie/Users_Guide.txt}{
              }{
        \logmessage{Directory 'Bibliographie' not found}
          }}}}}
\numberwithin{equation}{section}
\begin{document}

\title{The Levenberg-Marquardt Iteration for Numerical Inversion of the Power Density Operator}

\date{\today}

\author[1]{G.~Bal ({gb2030@columbia.edu})}

\author[2]{W.~Naetar ({wolf.naetar@univie.ac.at})}

\author[2,3]{O.~Scherzer ({otmar.scherzer@univie.ac.at})}

\author[4]{J.~Schotland ({schotland@umich.eu})}        

\affil[1]{\footnotesize Department of Applied Physics and Applied Mathematics, Columbia University \authorcr New York, NY 10027, USA \vspace{.5\baselineskip}}
\affil[2]{\footnotesize Computational Science Center, University of Vienna \authorcr Nordbergstr. 15, A-1090 Vienna, Austria \vspace{.5\baselineskip}}
\affil[3]{\footnotesize Radon Institute of Computational and Applied Mathematics, Austrian Academy of Sciences \authorcr Altenbergerstr. 69, A-4040 Linz, Austria \vspace{.5\baselineskip}}
\affil[4]{\footnotesize Department of Mathematics, University of Michigan \authorcr Ann Arbor, MI 48109, USA}

\maketitle

\section*{Abstract}
In this paper we develop a convergence analysis in an \emph{infinite} dimensional setting of the Levenberg-Marquardt iteration 
for the solution of a \emph{hybrid conductivity imaging} problem. The problem consists in determining the spatially varying conductivity $\sigma$ 
from a series of measurements of power densities for various voltage inductions. 
Although this problem has been very well studied in the literature, convergence and regularizing properties of iterative algorithms in an 
infinite dimensional setting are still rudimentary. We provide a partial result under the assumptions that the 
derivative of the operator, mapping conductivities to power densities, is injective and the data is noise-free. 
Moreover, we implemented the Levenberg-Marquardt algorithm and tested it on simulated data.

\paragraph*{Keywords.} Inverse problems, nonlinear ill-posed problems, iterative regularization, elliptic equations, hybrid imaging

\paragraph*{AMS subject classifications.} 35R30, 47J06, 35J47

\section{Introduction}
\label{sec:intro}

A common problem in \emph{hybrid imaging} consists in the determination of the spatially varying conductivity $\sigma > 0$ in a domain $\Omega \subset \Rm^n$ from $m$ 
measurements of power densities $\E_i(\sigma)=\sigma|\nabla u_i(\sigma)|^2$ ($i=1,\ldots,m$) inside $\Omega$ (resulting from $m$ different injected currents $f_i$).
That is, the potentials $u_i$ satisfy the elliptic equation

\begin{equation} \label{eq:calderon}
 \begin{aligned}
	\div (\sigma \nabla u_i) &= 0 \text{ in } \Omega\,, \\
	u_i &= f_i \text{ on } \partial\Omega.
 \end{aligned}
\end{equation}

This problem is relevant, for example, in \emph{Acousto-Electrical Tomography (AET)} \cite{ZhaWan04,AmmBonCapTanFin08,KucKun10} and \emph{Impedance-Acoustic Tomography (IAT)} \cite{GebSch08}.

We investigate the solution of this nonlinear inverse problem in an infinite-dimensional setting. We apply the \emph{Levenberg-Marquardt iteration} \cite{Han97}, a well-known iteration method which we recap in the next section. Using recent results about the linearized power density operator \cite{Bal12b}, we analyze local convergence conditions of the iteration method (for sufficiently smooth $\sigma$ and noise-free data $\E_i$) and provide a partial result. 

There are a number of theoretical results available on the problem of estimating $\sigma$ from power densities (and some additional boundary information). 
In a paper by Capdeboscq et al. \cite{CapFehGouKav09} it was shown that in $\mathbb{R}^2$, the conductivity $\sigma$ is uniquely determined by measurements 
$\F(\sigma)=(\sigma|\nabla u_1|^2,\sigma|\nabla u_2|^2,\sigma \nabla u_1 \nabla u_2)$ if

\begin{equation}
\label{eq:detcond}
	\det(\nabla u_1, \nabla u_2) \geq c > 0 \quad \text{in} \ \Omega.
\end{equation}

(certain Dirichlet boundary conditions are known to enforce this condition, see \cite{AleNes01}). Note that $\sigma \nabla u_1 \nabla u_2$ can easily be obtained from a third measurement of the power density using the polarization identity. \cite{BalBonMonTri11} extended this result to $\mathbb{R}^3$ and \cite{BalMon12} to arbitrary dimension (where the above determinant condition \eqref{eq:detcond} is much harder to fulfil) and additionally showed Lipschitz-stability of the reconstruction for sufficiently regular conductivities. Using the same interior measurements as above, Lipschitz-stability of the linearized problem $F'(\sigma)\colon L^2(\Omega) \to L^2(\Omega)$ (where $\sigma \in C^\infty(\Omega)$ and $\Omega \subset \mathbb{R}^2$) modulo the kernel of the linearized power density operator was shown by Kuchment and Steinhauer \cite{KucSte11}. What may be reconstructed in the setting of only {\em one} measurement $\sigma|\nabla u|^2$ is analyzed in \cite{BalAPDE12}.

Numerically, the problem has been treated in \cite{AmmBonCapTanFin08,GebSch08,CapFehGouKav09,KucKun11b,BalMon-2-12}. 

The paper is organized as follows: First, we introduce the \emph{Levenberg-Marquardt iteration} and its convergence conditions. Then, we linearize the power density operator $\E$ with Dirichlet boundary conditions in a suitable topology, analyse its stability and injectivity and discuss convergence of the iteration method. In the last two sections, we describe our numerical implementation of the Levenberg-Marquardt method and present numerical results.

\section{Iterative solution scheme}
\label{sec:itsol}

\subsection{The Levenberg-Marquardt iteration}
\label{sec:levmarq}

Let us denote by $\boldE^\delta = (\E_i^\delta)_{i=1,\ldots,m}$ the (noisy) measurements according to $m$ different initializations. We want find $\sigma$ such that $\boldE(\sigma) \approx \boldE^\delta$. To solve, we use an iterative regularization method in a Hilbert space setting. In \cite{KalNeuSch08,KalSchoSchu09,Kab11} several such methods for solving a general operator equation $F(\sigma) \approx y^\delta$, where $F\colon X \to Y$ and $X$, $Y$ are Hilbert spaces, are given. We use a Newton-type method, since these methods are usually converging faster than pure gradient descent methods.  Newton's method itself is defined by

\begin{equation*}
\label{eq:levmarq1}
	\begin{aligned}	
		\sigma_{k+1} = \sigma_k + F'(\sigma_k)^{-1}(y^\delta-F(\sigma_k))\;.
	\end{aligned}			
\end{equation*}

If the operator $F'(\sigma_k)$ is not left-invertible (this can be the case for $\boldE(\sigma)$ if there are not enough measurements, as we shall see in the next sections) one may use Tikhonov regularization:

\begin{equation*}
\label{eq:levmarq2}
	\sigma_{k+1} = \arg\min_\sigma \left\{\norm{y - F(\sigma_k) - F'(\sigma_k)(\sigma - \sigma_k)}_{Y}^2 + \alpha_k \norm{\sigma - \sigma_k}_{X}^2\right\}.
\end{equation*}

The resulting iteration method

\begin{equation}
\label{eq:levmarq3}
	\sigma_{k+1} = \sigma_k + ( F'(\sigma_k)^*F'(\sigma_k) + \alpha_k \Id )^{-1}F'(\sigma_k)^* (y^\delta-F(\sigma_k)),
\end{equation}

where $\alpha_k > 0$ are regularization parameters, is called the \emph{Levenberg-Marquardt method}.

Hanke \cite{Han97} suggests choosing $\alpha_k$ as the solutions of

\begin{equation}
\label{eq:hankecond}
	\norm{y^\delta - F(\sigma_k) - F'(\sigma_k)(\sigma_{k+1}(\alpha_k) - \sigma_k)}_Y = q \norm{y^\delta-F(\sigma_k))}_Y
\end{equation}

for some $0 < q < 1$. He shows that with the above choice of parameters and the initial value $\sigma_0$ sufficiently close to the desired solution, the Levenberg-Marquardt method converges monotonically to a solution $\sigma^\delta$ of $F(\sigma^\delta)=y^\delta$, provided the Fr\'{e}chet derivative $F'$ is uniformly bounded in a ball $B_\rho(\sigma^\delta)$ containing the initial value $\sigma_0$ and

\begin{equation}
\label{eq:conecond}
	\norm{F(\sigma) - F(\tilde{\sigma}) - F'(\sigma)(\sigma-\tilde{\sigma})}_Y \leq C \norm{\sigma-\tilde{\sigma}}_X \norm{F(\sigma)-F(\tilde{\sigma})}_Y
\end{equation}

for all $\sigma,\tilde{\sigma} \in B_\rho(\sigma^\delta)$.

\subsection{Linearization}
\label{sec:linear}

%In this section, we prove that the linearization given in \eqref{eq:syst1} is actually the Fr\'{e}chet derivative for appropriate spaces. We want $\E_i(\sigma)=\sigma|\nabla u_i(\sigma)|^2, i=1,\ldots,m$ to be in $L^2(\Omega)$ since the data might contain noise (for which differentiability cannot be assumed). 

%$\E(\sigma)$ is well-defined for $\sigma\colon \bar{\Omega} \to \mathbb{R}^+ \in L^\infty(\bar{\Omega})$ (from here on, we will always assume that $\sigma(x)$ is positive for all $x \in \Omega$). Indeed, let $L_\sigma u := \div(\sigma \nabla u)$ and $f \in H^{3/2}(\partial \Omega)$, then

%\begin{equation}
%	\label{eq:calderon2}
%	\begin{aligned}	
%		L_\sigma u &= 0 \quad \text{in} \ \Omega \\
%		u &= g \quad \text{in} \ \partial\Omega \\
%	\end{aligned}
%\end{equation}

We now prove Fr\'{e}chet differentiability of the power density operator $\E(\sigma)=\div (\sigma \nabla u)$ (where $u$ solves \eqref{eq:calderon} for some boundary condition $f=f_i$).

Let $\Omega \subset \Rm^n$. We assume that $\sigma$ is positive and in the space $H^{l}(\Omega)$ for $l>\frac n2$, which implies that $H^l(\Omega)$ is a Banach algebra with respect to pointwise multiplication \cite{Ada75}.  For $f \in H^{l +\frac 1 2}(\partial\Omega)$, \eqref{eq:calderon} then has a unique solution $u(\sigma) \in H^{l+1}(\Omega)$ \cite{Sol73}, so $\E\colon H^l(\Omega) \to H^l(\Omega)$.

Let $L_\sigma\colon H^l \to H^{l-2}, f \mapsto \div (\sigma \nabla f)$. By formal differentiation one can see that the directional derivative $u'(\sigma)\tau$ is given by

\begin{equation}
\label{eq:deriv_u}
	\begin{aligned}	
		L_\sigma u'(\sigma)\tau &= - \div(\tau \nabla u) \quad \text{in} \ \Omega \\
		u'(\sigma)\tau &= 0 \quad \text{on} \ \partial\Omega.
	\end{aligned}
\end{equation}

Obviously, $u'(\sigma)\colon H^l(\Omega) \to H^{l+1}(\Omega)$ is linear with respect to $\tau$. It is also bounded for the given norms considering

\begin{equation}
\label{eq:deriv_bounded}
	\norm{u'(\sigma)\tau}_{H^{l+1}(\Omega)} \leq C \norm{\div (\tau \nabla u)}_{H^{l-1}(\Omega)} \leq C \norm{\tau}_{H^l(\Omega)} \norm{u}_{H^{l+1}(\Omega)}.
\end{equation}

Here we used the regularity estimates in \cite{Sol73} for the first and the Banach algebra property for the last inequality. Since the constant in the regularity estimates is a bounded function of $\sigma$ if $\sigma$ is bounded from below (e.g., in a sufficiently small ball around the solution), $u'(\sigma)$ is actually uniformly bounded there. Now, let $R_u(\sigma,\tau)=u(\sigma + \tau) - u(\sigma) - u'(\sigma)\tau$ be the first order Taylor remainder of $u$. For small $\tau$ we get

\begin{equation*}
\label{eq:deriv_proof1}
		L_{\sigma+\tau} R_u(\sigma,\tau) = - L_\tau u'(\sigma)\tau \quad \text{in} \ \Omega.
\end{equation*}

As in \eqref{eq:deriv_bounded} we obtain

\begin{equation}
\label{eq:deriv_proof2}
	\begin{aligned}
		\norm{R_u(\sigma,\tau)}_{H^{l+1}(\Omega)} &\leq C \norm{L_\tau u'(\sigma)\tau}_{H^{l-1}(\Omega)} \\
		&\leq C \norm{\tau}_{H^l(\Omega)} \norm{u'(\sigma)\tau}_{H^{l+1}(\Omega)} \leq C \norm{\tau}_{H^l(\Omega)}^2
	\end{aligned}
\end{equation}

since $u'(\sigma)$ is bounded. This shows that $u'(\sigma)$ as defined in \eqref{eq:deriv_u} is actually the Fr\'{e}chet derivative of $u\colon H^l(\Omega) \to H^{l+1}(\Omega)$.

For the power density operator $\E\colon H^l(\Omega) \to H^l(\Omega), \ \sigma \mapsto \sigma |\nabla u(\sigma)|^2$, formal differentiation gives the directional derivative

\begin{equation*}
\label{eq:deriv_E}
	\E'(\sigma)\tau =  |\nabla u|^2 \tau + 2 \sigma \nabla u \nabla u'(\sigma)\tau,
\end{equation*}

Its uniform boundedness follows from the uniform boundedness of $u'(\sigma)$ and the Banach algebra property of $H^l(\Omega)$.

For the first order Taylor remainder $R_{\E}(\sigma,\tau)$, we get

\begin{equation*}
\label{eq:deriv_proof3}
	\begin{aligned}	
		R_{\E}(\sigma,\tau) &:= \E(\sigma+\tau)-\E(\sigma)-\E'(\sigma)\tau \\
		&= (\sigma+\tau) |\nabla(u + u'(\sigma)\tau + R_u(\sigma,\tau))|^2 - \sigma|\nabla u|^2 \\
		&- |\nabla u|^2 \tau - 2 \sigma \nabla u \nabla u'(\sigma)\tau \\
		&= (\sigma+\tau)\nabla( 2 u + 2 u'(\sigma)\tau + R_u(\sigma,\tau)) \nabla R_u(\sigma,\tau) \\
		&+ (\sigma+\tau)|\nabla u'(\sigma)\tau|^2 + 2 \tau \nabla u \nabla u'(\sigma)\tau.
	\end{aligned}	
\end{equation*}

Now, by \eqref{eq:deriv_proof2} and \eqref{eq:deriv_bounded}

\begin{equation}
\label{eq:deriv_proof4}
 \begin{aligned}
		\norm{R_{\E}(\sigma,\tau)}_{H^l(\Omega)} &\leq C (\norm{R_u(\sigma,\tau)}_{H^{l+1}(\Omega)} + \norm{\tau}_{H^l(\Omega)}^2 + \norm{\tau}_{H^l(\Omega)}^2) \\
		&\leq C \norm{\tau}_{H^l(\Omega)}^2
 \end{aligned}
\end{equation}

so $\norm{R_{\E}(\sigma,\tau)}_{H^l(\Omega)} = o(\norm{\tau}_{H^l(\Omega)})$, showing that $\E'(\sigma)$ is the Fr\'{e}chet derivative of $\E\colon H^l(\Omega) \to H^l(\Omega)$ at $\sigma$.

\subsection{Stability and injectivity}
\label{sec:stinj}

In this section, we recap stability and injectivity of the linearized operator, following the treatment in \cite{Bal12b}. 

We may write the linearization of $\boldE\colon H^l(\Omega) \to H^l(\Omega;\Rm^m), \sigma \mapsto (\E_i)_{i=1,\ldots,m}$ as the following system:

\begin{equation} \label{eq:syst1}
 \begin{aligned}
	\nabla\cdot\delta\sigma \nabla u_i + \nabla\cdot \sigma \nabla \delta u_i &= 0 \\
	\delta\sigma|\nabla u_i|^2 + 2\sigma\nabla u_i\cdot\nabla\delta u_i &= \delta {\E}_i
 \end{aligned}
\end{equation}

The last line provides the Fr\'echet derivative $\delta{\E}_i=d{\E}_i \cdot \delta\sigma$ seeing $\delta u_i=du_i\cdot\delta\sigma$ as a function of $\delta\sigma$ given by solving the first line with Dirichlet conditions. 

\paragraph*{Ellipticity.}
Let us assume as before that $\sigma$ and $\nabla u_i$ are in the space $H^{l}(\Omega)$ for $l>\frac n2$. Then Dirichlet conditions for $\delta u_i$ render the above system elliptic provided that $q_i(x,\xi)=0$ for all $1\leq j\leq m$ implies $\xi=0$, where 

\begin{equation*}
  q_i(x,\xi) = 2(\hat F_i\cdot \xi)^2 - |\xi|^2,\qquad \hat F_i =\dfrac{\nabla u_i}{|\nabla u_i|}.
\end{equation*}

This theory is independent of dimension and is based on the theory of elliptic redundant systems of equations with boundary conditions satisfying the Lopatinskii conditions; see \cite{Bal12b,Sol73}. Similar ellipticity conditions are found in \cite{KucSte11}. From such a theory, we obtain for $s=l$ with $l>\frac n2$, and coefficients $\sigma\in H^{s}(\Omega)$ and $u_i\in H^{s+1}(\Omega)$,  the following stability estimate

\begin{equation} \label{eq:stab1}
 \begin{aligned}
  \|\delta \sigma - \delta \tilde\sigma\|_{H^s(\Omega)} &+ \| \delta \mathbf{u}- \delta \tilde{\mathbf{u}}\|_{H^{s+1}(\Omega;\Rm^m)} \\
  &\leq C \| \delta \boldE- \delta \tilde{\boldE}\|_{H^s(\Omega;\Rm^m)} + C_2 \|\delta \mathbf{u}-\delta \tilde{\mathbf{u}} \|_{L^2(\Omega;\Rm^m)}.
  \end{aligned}
\end{equation}

Here, $\delta\boldE$ is the collection of terms $\delta\E_i$ and $\delta \mathbf{u}$ the collection of $\delta u_i$. The presence of the term with $C_2$ indicates the fact that the systems may not be invertible. Rather, it is Fredholm, and hence invertible up to a finite dimensional kernel of (smooth) functions. This is reminiscent of the behavior of $\Delta+q(x)$, which is invertible up to a finite dimensional kernel, although it is invertible (with Dirichlet conditions) for most values of $q(x)$.

\paragraph*{Injectivity of modified systems.}
Whether we can choose $C_2=0$ above is not known in general. Injectivity results are known in the presence of more measurements than necessary for ellipticity; see, e.g., \cite{KucKun11b} in the linearized setting and \cite{BalBonMonTri11,BalMon12} in the nonlinear setting.
%This is the case in the vicinity of $\sigma$ a constant with appropriate choices of $u_i$; see \cite{KK-AET-11}. NO: PROOF REQUIRES MORE MEASUREMENTS.
 
In the setting of two measurements in dimension $n=2$ and three measurements in dimension $n=3$, results of injectivity were obtained in \cite{Bal12b} for modified systems of equations. Let denote by $v=(\delta\sigma, \{\delta u_i\})$ and recast the system \eqref{eq:syst1} as 

\begin{equation*}
 {\mathcal A} v ={\mathcal S}.
\end{equation*}

Then assuming that full Cauchy data are known for $v$ on $\partial \Omega$, we may replace the above system by

\begin{equation}
\label{eq:syst2}
{\mathcal A}^t{\mathcal A} v ={\mathcal A}^t{\mathcal S}\quad\text{ in }\Omega, \qquad v= \partial_\nu v=0 \text{ on } \partial \Omega.
\end{equation}

This system was proved in \cite{Bal12b} to admit a unique solution on sufficiently small domains; in other words, \eqref{eq:stab1} holds with $C_2=0$.

A different modification of \eqref{eq:syst1} also leading to a fourth-order system similar to \eqref{eq:syst2} is also proved to be injective; see \cite[Equation (35)]{Bal12b}.

These systems are difficult to implement numerically. We therefore wish to solve the system as written in \eqref{eq:syst1}. Knowing that very similar systems are indeed shown to be injective in the aforementioned work, it seems reasonable to assume that $C_2=0$ in \eqref{eq:stab1} for $m$ large enough, i.e., that the linearized operator ${\mathcal A}$ is left-invertible. With this assumption, we get that $d{\boldE}$ is an operator from $H^s(\Omega)$ to $H^s(\Omega;\Rm^m)$ with a bounded left-inverse when $s>\frac n2$.

We then obtain an optimal stability estimate of the form 

\begin{equation} \label{eq:linear_stability}
  C_s^{-1} \|\delta\sigma-\delta\tilde\sigma\|_{H^s(\Omega)} \leq  \|\delta{\boldE}-\delta \tilde{\boldE}\|_{H^s(\Omega;\Rm^m)} \leq C_s \|\delta\sigma-\delta\tilde\sigma\|_{H^s(\Omega)}.
\end{equation}

\subsection{Convergence analysis}
\label{sec:convergence}

In this section, we analyse whether the convergence conditions given in section \ref{sec:levmarq} can be shown to hold for the power density operator $\boldE\colon H^{l}(\Omega) \to H^l(\Omega;\mathbb{R}^2)$ with a suitable number of measurements $m$. Note that we require noise-free data here (the noise cannot be assumed fulfil differentiability constraints).

To get local convergence of the Levenberg-Marquardt iteration, it suffices to show that \eqref{eq:conecond} holds (we already saw that the operator is uniformly bounded close to a solution). Given Lipschitz-type stability, i.e.,

\begin{equation}
\label{eq:conv_stability_ideal}
	\norm{\sigma_1 - \sigma_2}_{H^l(\Omega)} \leq C' \norm{\boldE(\sigma_1)-\boldE(\sigma_2)}_{H^l(\Omega;\mathbb{R}^m)},
\end{equation}

the estimate \eqref{eq:conecond} can easily be obtained from \eqref{eq:deriv_proof4}.

Using the stability estimates \eqref{eq:linear_stability} for the linearized operator, we can get such an estimate locally, since for $\sigma_1,\sigma_2 \in B_\rho(\sigma^\delta)$ (a ball with sufficiently small diameter in the $H^l(\Omega)$-norm) we have

\begin{equation*}
\label{eq:conv_stability_proof1}
	\begin{aligned}	
		C' \norm{\sigma_2-\sigma_1}_{H^l(\Omega)} &\leq C_s^{-1} \norm{\sigma_2-\sigma_1}_{H^l(\Omega)} - C \norm{\sigma_2-\sigma_1}_{H^l(\Omega)}^2 \\
		&\leq \norm{\boldE'(\sigma_1)(\sigma_1 - \sigma_2)}_{H^l(\Omega;\mathbb{R}^m)} - C \norm{\sigma_2-\sigma_1}_{H^l(\Omega)}^2 \\
		&\leq \norm{\boldE(\sigma_2)-\boldE(\sigma_1)}_{H^l(\Omega;\mathbb{R}^m)}
	\end{aligned}	
\end{equation*}

by using norm equivalence for the second inequality and 

\begin{equation*}
\label{eq:conv_stability_proof2}
	\norm{\boldE'(\sigma_1)(\sigma_2-\sigma_1)}_{H^l(\Omega;\mathbb{R}^m)} - \norm{\boldE(\sigma_2)-\boldE(\sigma_1)}_{H^l(\Omega;\mathbb{R}^m)} 
	\leq C \norm{\sigma_2-\sigma_1}_{H^l(\Omega)}^2,
\end{equation*}

which is obtained by applying the reverse triangle inequality to \eqref{eq:deriv_proof4}, for the third estimate.

From \eqref{eq:deriv_proof4} we finally get

\begin{equation}
\label{eq:conv_condition}
	\norm{R_{\boldE}(\sigma,\tau)}_{H^l(\Omega;\mathbb{R}^m)} \leq C \norm{\tau}_{H^l(\Omega)} \norm{\boldE(\sigma+\tau)-\boldE(\sigma)}_{H^l(\Omega;\mathbb{R}^m)}.
\end{equation}

for sufficiently small $\tau$. Applying the general results in \cite{Han97} now gives us local convergence (in the noise-free case) of the Levenberg-Marquardt iteration applied to $\boldE\colon H^l(\Omega) \to H^l(\Omega;\mathbb{R}^m)$.

Note that for the stability estimate \eqref{eq:linear_stability} we assumed injectivity of \eqref{eq:syst1} (and thus $\E'(\sigma)$). For one measurement, this may not hold, for example, it is easily seen that if $\Omega$ is the unit circle in $\Rm^2$, the boundary voltage $f(x,y)=ax+by$ (where $\vectwo{a}{b}$ is a unit vector) and $\sigma=1$, the function $f$ (with corresponding linearized potential $(u'(\sigma)f)(x,y)=-\frac{1}{4}(x^2+y^2)+\frac{1}{4}$) is in $\operatorname{Ker}\E'(\sigma)$; see also \cite{BalAPDE12}.

Generically, one may however expect injectivity to hold if enough data are available. 

\section{Numerical solution}
\label{sec:numerical}

The following sections contain some numerical work to demonstrate the feasibility of our ideas. For simplicity, we work in two dimensions, so we take $\Omega \subset \Rm^2$. Furthermore, we now allow the presence of noise. As noise cannot be assumed to be differentiable, we have to consider (in contrast to the previous section) the power density operator to map into $L^2(\Omega)$. Since $l=2$ is the smallest integer for which $H^l(\Omega)$ is a Banach algebra, we take $\E\colon H^2(\Omega) \to L^2(\Omega)$.

\subsection{Calculation of the adjoint}
\label{sec:adjoints}

For the iterative algorithm, we require an expression for the adjoint of $\E'(\sigma)$, which we now derive.

First, we consider the case $m=1$ of a single measurement. Larger values of $m$ are treated with an additional summation. Since $\E'(\sigma)\colon H^2(\Omega) \to L^2(\Omega)$ is a bounded operator, we know that $\E'(\sigma)^*\colon L^2(\Omega) \to H^2(\Omega)$ is bounded as well (see \cite{Wei80}).

First note that the adjoint operator of $\E'(\sigma)\colon H^2(\Omega) \to L^2(\Omega)$ can be written as

\begin{equation*}
\label{eq:adjoints_proof1}
	\begin{aligned}	
		&\E'(\sigma)^*\colon L^2(\Omega) \to H^2(\Omega) \\
		&\E'(\sigma)^* z = i^* \tilde{\E'}(\sigma)^* z,
	\end{aligned}	
\end{equation*}

where $\tilde{\E}'(\sigma)^*$ is the $L^2(\Omega)$-adjoint of $\E(\sigma)$ and $i^*$ the adjoint of the embedding $i\colon H^2(\Omega) \to L^2(\Omega)$.

We start by finding an expression for $\tilde{\E}'(\sigma)^*$ on $H^2(\Omega) \subset L^2(\Omega)$. Let $V=V(\sigma,u(\sigma))$ with $V\colon H^2(\Omega) \to H^2_0(\Omega)$ be the linear operator defined by 

\begin{equation*}
\label{eq:adjoints_proof2}
	\begin{aligned}
		L_\sigma V z  &= -\div(z \sigma \nabla u(\sigma) ) \quad \text{in} \ \Omega \\
		V z &= 0 \quad \text{in} \ \partial\Omega .\\
	\end{aligned}	
\end{equation*}

Note that $z\sigma$ is also in $H^2(\Omega)$, so $V$ is mapped into $L^2(\Omega)$.

For $z \in H^2$ we now get that for the second summand of $\E'(\sigma)$

\begin{equation*}
\label{eq:adjoints_proof3}
	\begin{aligned}	
		\iprod{\sigma \nabla u \nabla u'(\sigma) \tau,z}_{L^2(\Omega)} &= -\iprod{\sigma \nabla u'(\sigma) \tau,\nabla V z}_{L^2(\Omega;\mathbb{R}^2)} \\
		&= \iprod{\tau,\nabla u \nabla V z}_{L^2(\Omega)}
	\end{aligned}			
\end{equation*}

by partial integration and the definitions of $V$ and $u'(\sigma)$. Since the first summand is self-adjoint, we conclude

\begin{equation}
\label{eq:adjoints_proof4}
	\begin{aligned}	
		&\E'(\sigma)^*\colon H^2(\Omega) \to H^2(\Omega) \\
		&\E'(\sigma)^* z = i^* ( |\nabla u(\sigma)|^2 z + 2 \nabla u(\sigma) \nabla V z )\;.
	\end{aligned}			
\end{equation}

As $\E'(\sigma)^*$ is bounded on $L^2(\Omega)$, \eqref{eq:adjoints_proof4} can be continuously extended to $w \in L^2(\Omega)$, e.g., by taking $\E'(\sigma)^*w:=\lim_{\epsilon \to 0} \E'(\sigma)^*\Phi_\epsilon w$ where $\Phi_\epsilon w = \phi_\epsilon \ast w$ is the mollification operator corresponding to some mollifier $\phi$.

For $\boldE\colon H^2(\Omega) \to L^2(\Omega;\mathbb{R}^m), \ \sigma \mapsto (\E_i)_{i=1,\ldots,m}$ we similarly get (with $V_j=V(\sigma,u_j(\sigma))$),

\begin{equation}
\label{eq:adjoints_vectorial}
	\begin{aligned}	
		&\boldE'(\sigma)^*\colon H^2(\Omega;\mathbb{R}^2) \to H^2(\Omega) \\
		&\boldE'(\sigma)^* z=i^* \sum_{j=1}^m |\nabla u_j(\sigma)|^2 z_j + 2 \nabla u_j(\sigma) \cdot \nabla V_j w_j \;.
	\end{aligned}			
\end{equation}

$\boldE'$ can also be continuously extended to $L^2(\Omega;\mathbb{R}^2)$ by mollifying and taking the limit.

To calculate $i^*\colon L^2(\Omega) \to H^2(\Omega)$ we make the additional assumption that $\sigma \in H^2_N(\Omega) = \{ x \in H^2(\Omega) \big| \dn{x} = 0 \}$, a closed subspace of $H^2(\Omega)$ since $(\frac{\partial}{\partial \nu} \circ \operatorname{trace})\colon H^2(\Omega) \to L^2(\partial \Omega)$ is bounded.

With the inner product

\begin{equation*}
\label{eq:adjoints_proof5}
	\iprod{x,y}_{H^2_N(\Omega)}=\iprod{x,y}_{L^2(\Omega)} + \beta^2 \iprod{\Delta x,\Delta y}_{L^2(\Omega)}
\end{equation*}

for some $\beta > 0$, $i^* y$ is the solution of the Neumann problem

\begin{equation}
\label{eq:adjoints_proof6}
	\begin{aligned}	
		(\Id + \beta^2 \Delta \Delta)i^* y &= y \\		
		\dn{\Delta i^*y} &= 0 \text{ on } \partial\Omega \\
		\dn{i^*y} &= 0 \text{ on } \partial\Omega.
	\end{aligned}			
\end{equation}

since for all $x \in H^2_N(\Omega)$, $y \in L^2(\Omega)$ we get that

\begin{equation*}
\label{eq:adjoints_proof7}
	\begin{aligned}	
		\iprod{x,i^*y}_{H^2_N(\Omega)} &= \iprod{x,i^*y}_{L^2(\Omega)} + \beta^2 \iprod{\Delta x,\Delta i^*y}_{L^2(\Omega)} \\
		&= \iprod{x,(\Id + \beta^2 \Delta \Delta) i^* y}_{L^2(\Omega)} + \beta^2 \int_{\partial\Omega} \dn{x} \Delta i^*y - \dn{\Delta i^*y} x \ dS \\
		&= \iprod{x,y}_{L^2(\Omega)} = \iprod{ix,y}_{L^2(\Omega)}.
	\end{aligned}			
\end{equation*}

The fourth order PDE \eqref{eq:adjoints_proof6} is equivalent to the system of second order equations

\begin{equation}
\label{eq:adjoints_proof8}
	\begin{aligned}	
		i^* y + \beta^2 \Delta z &= y \\
		\Delta i^* y - z &= 0 \\
		\dn{z} &= 0 \text{ on } \partial\Omega \\
		\dn{i^*y} &= 0 \text{ on } \partial\Omega.
	\end{aligned}			
\end{equation}

For a cruder approximation (which should give a faster solution), we could use $H^1$ or $L^2$ adjoints (where we have $i=i^*=\operatorname{Id}$). For the inner product $\iprod{x,y}_{H^1(\Omega)}=\iprod{x,y}_{L^2(\Omega)} + \beta \iprod{\nabla x,\nabla y}_{L^2(\Omega)}$, the adjoint of $\tilde{i}\colon H^1(\Omega) \to L^2(\Omega;\mathbb{R}^2)$ (by a proof similar to the above) maps $y \in L^2(\Omega)$ to the solution of

\begin{equation}
\label{eq:adjoints_proof9}
	\begin{aligned}	
		(\Id - \beta^2 \Delta) \tilde{i}^* y &= y \\		
		\dn{\tilde{i}^*y} &= 0 \text{ on } \partial\Omega.
	\end{aligned}			
\end{equation}

\subsection{Implementation}
\label{sec:implem}

For $F=\E$ with a given (single) measurement $\E^\delta$, the equation \eqref{eq:levmarq3} for the $k$-th Levenberg-Marquardt step $\tau_k$ reads

\begin{equation}
\label{eq:num_methods1}
	(\lim_{\epsilon \to 0} i^*M_\epsilon + \alpha_k \Id)\tau_k = \E'(\sigma_k)^*(\E^\delta - \E(\sigma_k)) \\
\end{equation}

where, with $u=u(\sigma_k)$, $\alpha_k=\alpha$, $\tau=\tau_k$ and $\sigma=\sigma_k$,

\begin{equation*}
\label{eq:num_methods2}
	\begin{aligned}	
		M_\epsilon \tau &= |\nabla u|^2 \Phi_\epsilon (|\nabla u|^2 \tau + 2 \sigma \nabla u \nabla u'(\sigma)\tau) \\
		&+ 2 \nabla u \nabla V \Phi_\epsilon (\tau |\nabla u|^2) + 4 \nabla u \nabla V\Phi_\epsilon (\sigma \nabla u \nabla u'(\sigma)\tau).
	\end{aligned}			
\end{equation*}

Now, let $\epsilon$ be fixed. By introducing auxiliary variables $z^1,z^2,z^3$ and $y^1,y^2,y^3$ (for expressions which are a solution of a PDE) and using \eqref{eq:adjoints_proof8}, equation \eqref{eq:num_methods1} can be written as

\small
\begin{equation}
\label{eq:num_methods3}
	\begin{alignedat}{2}
		&z^3 + \alpha \tau_\epsilon = \E'(\sigma)^*(\E^\delta - \E(\sigma)) \\
		\\
		&\Delta z^3 - z^2 = 0, &\hspace{-2.4 cm} \dnflat{z^3} = 0 \text{ on } \partial\Omega \\ 				
		&\beta^2 \Delta z^2 + z^3 - z^1 = 0, &\hspace{-2.4 cm} \dnflat{z^2} = 0 \text{ on } \partial\Omega \\		
		\\
		&|\nabla u|^2 \Phi_\epsilon ( |\nabla u|^2 \tau_\epsilon + 2 \sigma \nabla u \nabla y^1) + 2 \nabla u \nabla y^2 + 4 \nabla u \nabla y^3 - z^1 = 0 \\		
		\\
 		&\div(\Phi_\epsilon (\sigma \nabla u \nabla y^1) \sigma \nabla u) + L_{\sigma} y^3 = 0, &\hspace{-2.7 cm}y^3 = 0 \text{ on } \partial\Omega \\		
		&\div(\Phi_\epsilon (\tau_\epsilon |\nabla u|^2) \sigma \nabla u) + L_{\sigma} y^2 = 0, &\hspace{-2.7 cm}y^2 = 0 \text{ on } \partial\Omega \\		
		&\div(\tau_\epsilon \nabla u) + L_{\sigma} y^1 = 0, &\hspace{-2.7 cm}y^1 = 0 \text{ on } \partial\Omega \\		
	\end{alignedat}		
\end{equation}

\normalsize
We then obtain $\tau_k$ in \eqref{eq:num_methods1} from $\tau_\epsilon$ by letting $\epsilon \to 0$ (for numerical purposes the continuous extension, i.e., the mollification and the limit can be ignored, as it leaves the discretization unchanged in the limit $\epsilon \to 0$). Since its discretization does not contain inverse matrices (and thus can be solved using sparse matrix operations only), the system \eqref{eq:num_methods3} is more amenable to numerical treatment than discretizing \eqref{eq:num_methods1} directly. When using $L^2$ or $H^1$ adjoints, similar (smaller) systems can be generated.

For multiple measurements $\boldE$, using the exact same approach leads to the system

\small
\begin{equation}
\label{eq:num_methods4}
	\begin{alignedat}{2}
		&z^3 + \alpha \tau_\epsilon = \boldE'(\sigma)^*(\boldE^\delta - \boldE(\sigma)) \\
		\\
		&\Delta z^3 - z^2 = 0, &\hspace{-3.8 cm} \dnflat{z^3} = 0 \text{ on } \partial\Omega \\ 				
		&\beta^2 \Delta z^2 + z^3 - z^1 = 0, &\hspace{-3.8 cm} \dnflat{z^2} = 0 \text{ on } \partial\Omega \\		
		\\
		&\sum_{j=1}^m |\nabla u_j|^2 \Phi_\epsilon ( |\nabla u_j|^2 \tau_\epsilon + 2 \sigma \nabla u_j \nabla y^1_j) + 2 \nabla u_j \nabla y^2_j + 4 \nabla u_j \nabla y^3_j - z^1 = 0 \\		
		\\
 		&\div(\Phi_\epsilon (\sigma \nabla u_j \nabla y^1_j) \sigma \nabla u_j) + L_{\sigma} y^3_j = 0,  &\hspace{-2.2 cm}y^3_j = 0 \text{ on } \partial\Omega \\
		&\div(\Phi_\epsilon (\tau_\epsilon |\nabla u_j|^2) \sigma \nabla u_j) + L_{\sigma} y^2_j = 0,  &\hspace{-2.2 cm}y^2_j = 0 \text{ on } \partial\Omega \\
		&\div(\tau_\epsilon \nabla u_j) + L_{\sigma} y^1_j = 0,  &\hspace{-2.2 cm}y^1_j = 0 \text{ on } \partial\Omega \\
	\end{alignedat}		
\end{equation}

\normalsize
so 3 equations (and auxiliary variables) have to be added per additional measurement. We again obtain $\tau_k$ from $\tau_\epsilon$ by letting $\epsilon \to 0$.

\section{Results}
\label{sec:results}

In our numerical solution we directly implemented the systems \eqref{eq:num_methods3} and \eqref{eq:num_methods4} to obtain Levenberg-Marquardt steps. For discretization of the partial differential equations we used a self-written linear finite elements framework in \textit{MATLAB}. The triangular mesh (with 42849 nodes and 85007 elements) was created using \textit{DistMesh}. \cite{PerStr04} All calculations were done on a workstation computer.

To test the reconstruction algorithm, we generated simulated data $\boldE^\delta=(\E_1^\delta,\E_2^\delta,\E_3^\delta)$ on the unit circle with the boundary conditions $f_1(x,y)=x$, $f_2(x,y)=y$ and $f_3(x,y)=\frac{1}{\sqrt{2}}(x-y)$ (and added Gaussian noise with standard deviation $1$ to avoid an inverse crime). 

\FloatBarrier
\begin{figure}[htp]
	\centering

		\subfloat[Conductivity $\sigma$]{\includegraphics[width=0.4\textwidth]{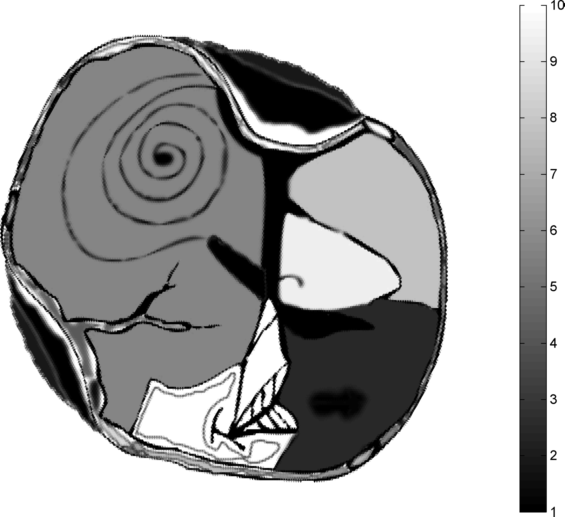}} \hspace{0.1\textwidth}
  		\subfloat[Power density $\E_1^\delta$]{\includegraphics[width=0.4\textwidth]{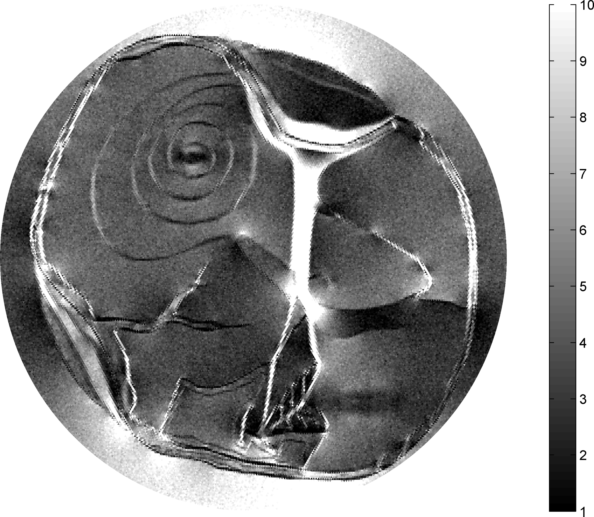}} \\
  		\subfloat[Power density $\E_2^\delta$]{\includegraphics[width=0.4\textwidth]{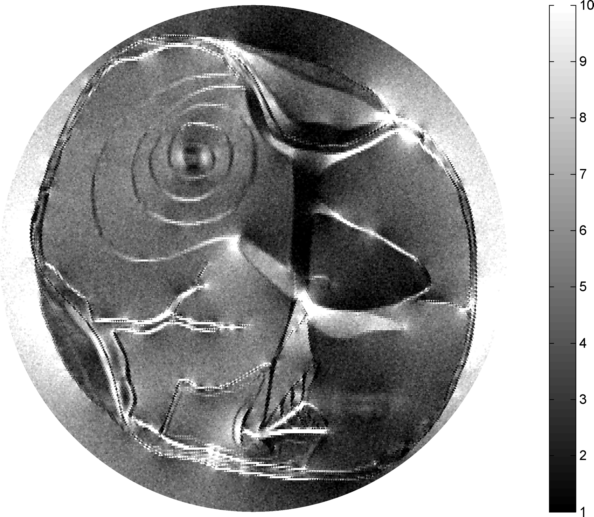}} \hspace{0.1\textwidth}
  		\subfloat[Power density $\E_3^\delta$]{\includegraphics[width=0.4\textwidth]{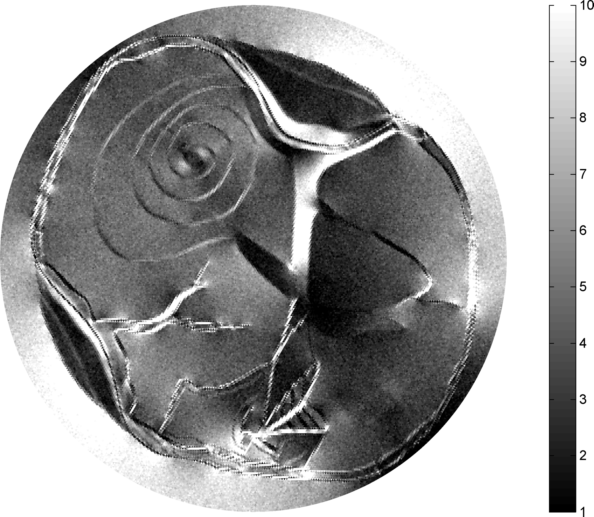}}

	\caption[Conductivity and power density]
	{
	  Conductivity and simulated measurements with boundary conditions $f_1(x,y)=x$, $f_2(x,y)=y$ and $f_3(x,y)=\frac{1}{\sqrt{2}}(x-y)$ and Gaussian noise. The color axis was manually set to 1-10, the range of $\sigma$.
	}	
	\label{fig:image_cond_pd}	
\end{figure}
\FloatBarrier

For reconstruction, we set $\sigma_0 \equiv 1$, $\beta=10^{-3}$ and chose $\alpha_k = \frac{1}{2^k}$ (with a minimum of $10^{-8}$) a priori. Choosing $\alpha_k$ according to the criterion \eqref{eq:hankecond} would be possible (e.g., by reducing $\alpha_k$ until the condition is met), but numerically expensive. 

Note that for two or more measurements, regularization is technically is not necessary since the operator can be assumed to be invertible (the linearized and discretized equations were solvable with $\alpha_k=0$ in all instances we tested), but advantageous for the iteration scheme since it helps ensure that iterates remain in a trust region and thus feasible (i.e., positive). Additionally, we enforce a minimal conductivity of $10^{-12}$ on the iterates. 

We used the standard \textit{MATLAB} sparse solver \textit{mldivide} to solve the linearized and discretized equations.

\begin{figure}[htp]
	\centering
		
  		\subfloat[$L^2$ adjoints, data $\E_2^\delta$]{\includegraphics[width=0.32\textwidth]{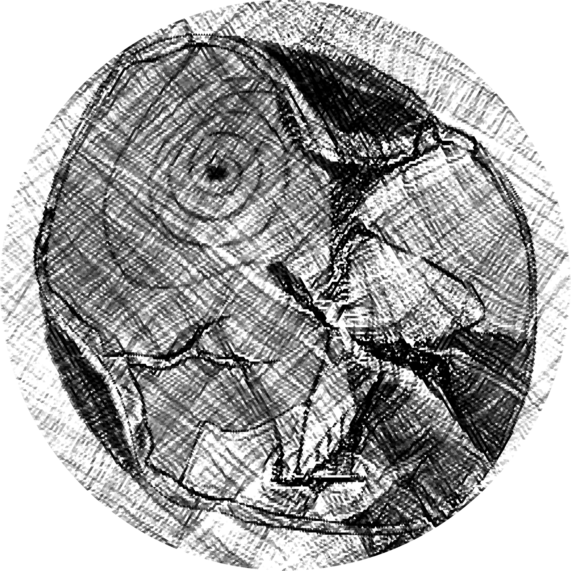}} \hspace{0.01\textwidth}
  		\subfloat[$L^2$ adjoints, data $(\E_1^\delta,\E_3^\delta)$]{\includegraphics[width=0.32\textwidth]{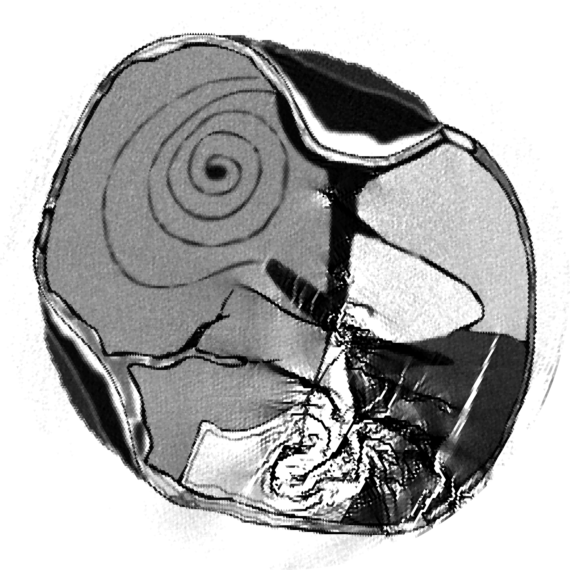}} \hspace{0.01\textwidth}
  		\subfloat[$L^2$ adjoints, data $(\E_1^\delta,\E_2^\delta,\E_3^\delta)$]{\includegraphics[width=0.32\textwidth]{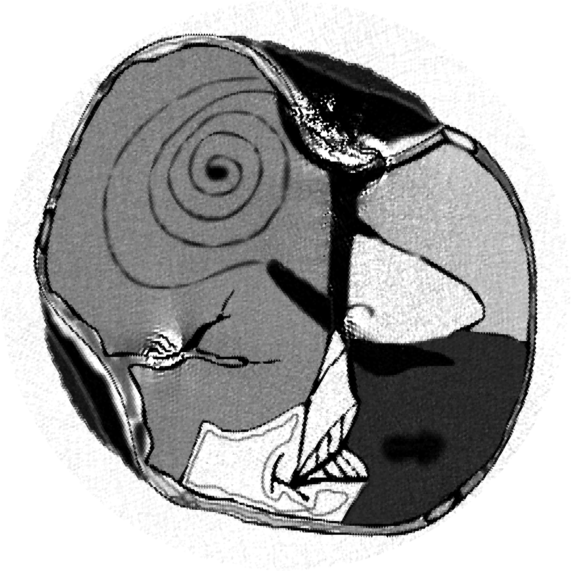}} \\
  		\subfloat[$H^1$ adjoints, data $\E_2^\delta$]{\includegraphics[width=0.32\textwidth]{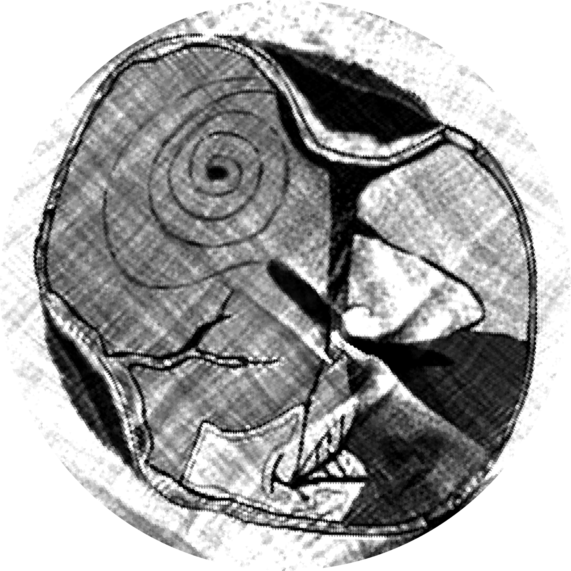}} \hspace{0.01\textwidth}
  		\subfloat[$H^1$ adjoints, data $(\E_1^\delta,\E_3^\delta)$]{\includegraphics[width=0.32\textwidth]{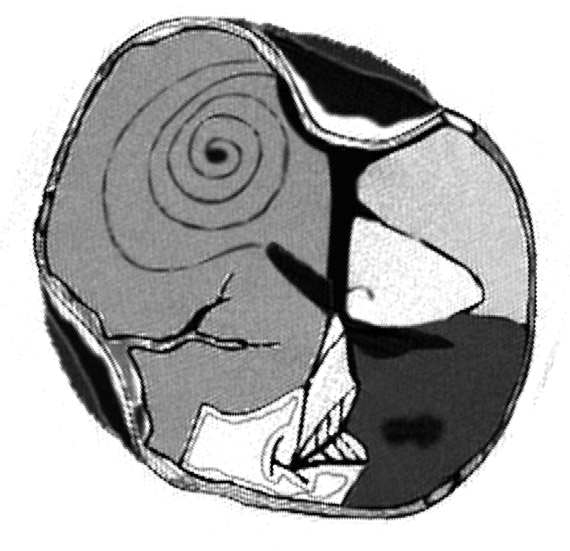}} \hspace{0.01\textwidth}
  		\subfloat[$H^1$ adjoints, data $(\E_1^\delta,\E_2^\delta,\E_3^\delta)$]{\includegraphics[width=0.32\textwidth]{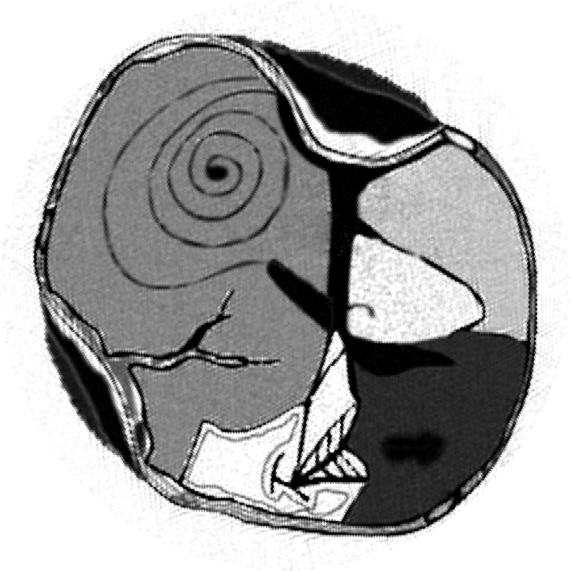}} \\
  		\subfloat[$H^2$ adjoints, data $\E_2^\delta$]{\includegraphics[width=0.32\textwidth]{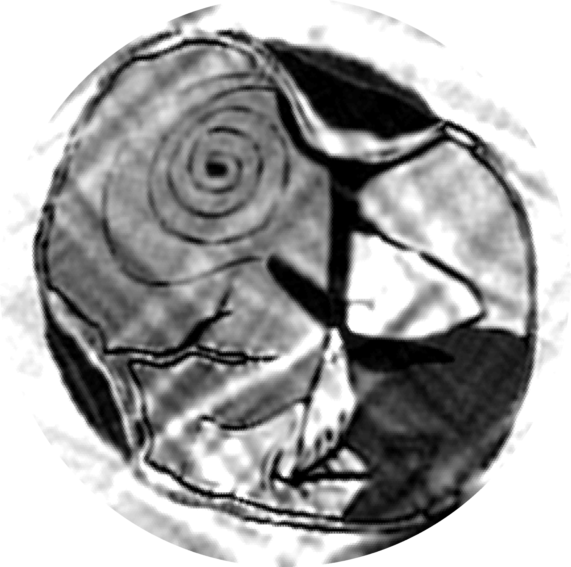}} \hspace{0.01\textwidth}
  		\subfloat[$H^2$ adjoints, data $(\E_1^\delta,\E_3^\delta)$]{\includegraphics[width=0.32\textwidth]{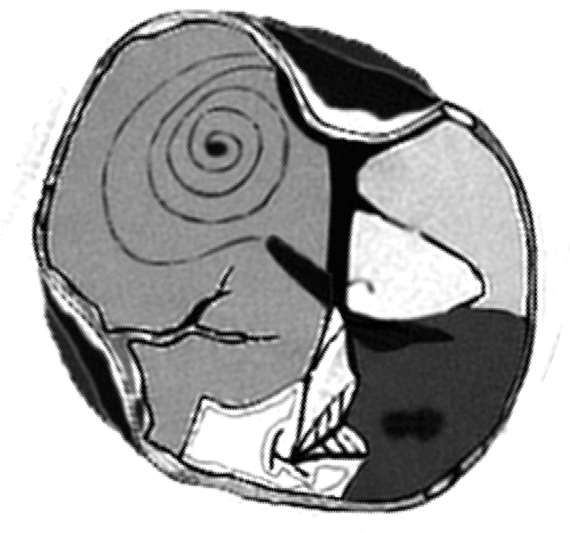}} \hspace{0.01\textwidth}
  		\subfloat[$H^2$ adjoints, data $(\E_1^\delta,\E_2^\delta,\E_3^\delta)$]{\includegraphics[width=0.32\textwidth]{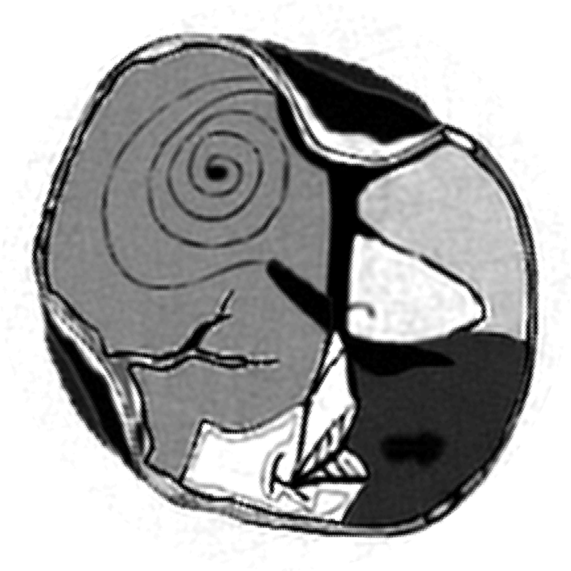}}

	\caption[Reconstructed conductivities]
	{
	  Reconstructed conductivities after 15 iterations of the Levenberg-Marquardt algorithm using different adjoints and one, two or three measurements. The color axis was manually set to 1-10, the range of the original conductivity to ease comparison.
	}	
	\label{fig:image_results}	
\end{figure}

\begin{figure}[htp]
	\centering
		
  		\subfloat[$L^2$ adjoints, data $\E_2^\delta$]{\includegraphics[width=0.32\textwidth]{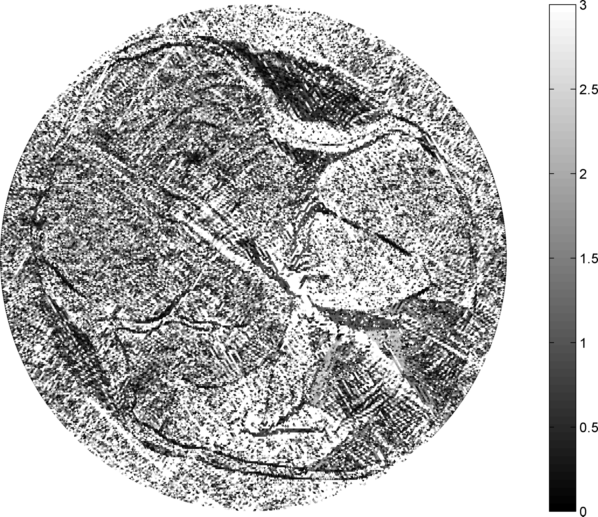}}  \hspace{0.01\textwidth}
  		\subfloat[$L^2$ adjoints, data $(\E_1^\delta,\E_3^\delta)$]{\includegraphics[width=0.32\textwidth]{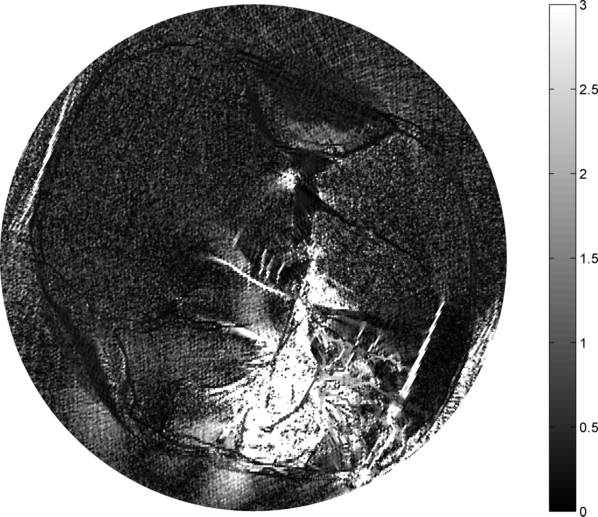}} \hspace{0.01\textwidth}
  		\subfloat[$L^2$ adjoints, data $(\E_1^\delta,\E_2^\delta,\E_3^\delta)$]{\includegraphics[width=0.32\textwidth]{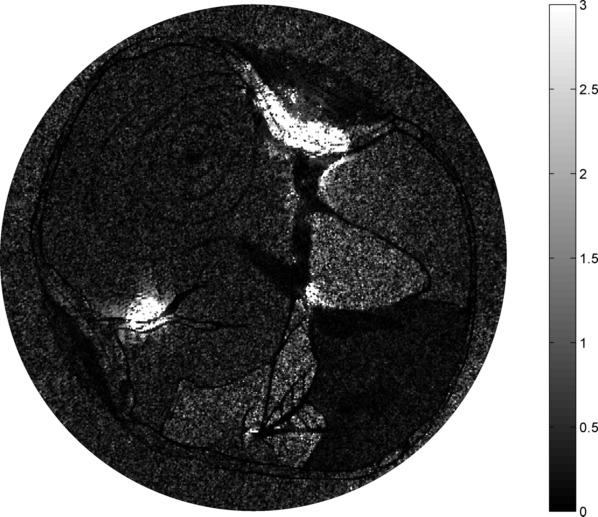}} \\
  		\subfloat[$H^1$ adjoints, data $\E_2^\delta$]{\includegraphics[width=0.32\textwidth]{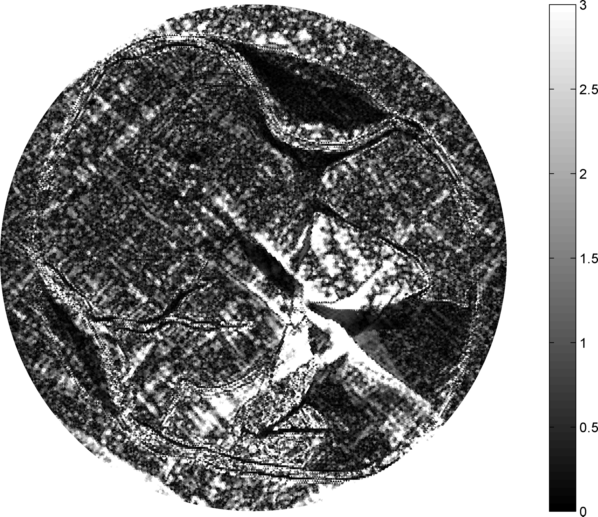}} \hspace{0.01\textwidth}
  		\subfloat[$H^1$ adjoints, data $(\E_1^\delta,\E_3^\delta)$]{\includegraphics[width=0.32\textwidth]{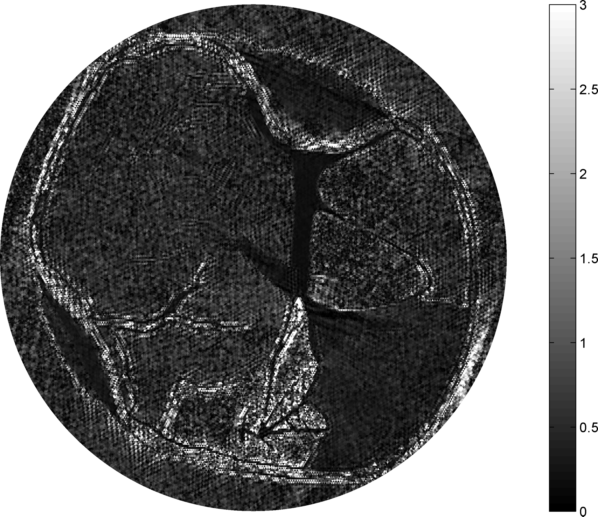}} \hspace{0.01\textwidth}
  		\subfloat[$H^1$ adjoints, data $(\E_1^\delta,\E_2^\delta,\E_3^\delta)$]{\includegraphics[width=0.32\textwidth]{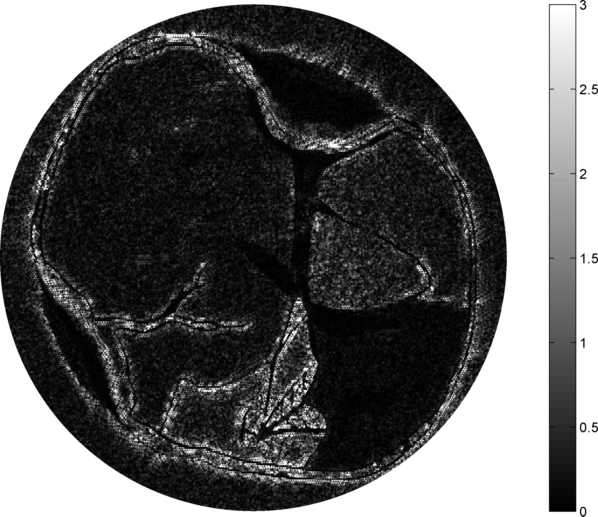}} \\
  		\subfloat[$H^2$ adjoints, data $\E_2^\delta$]{\includegraphics[width=0.32\textwidth]{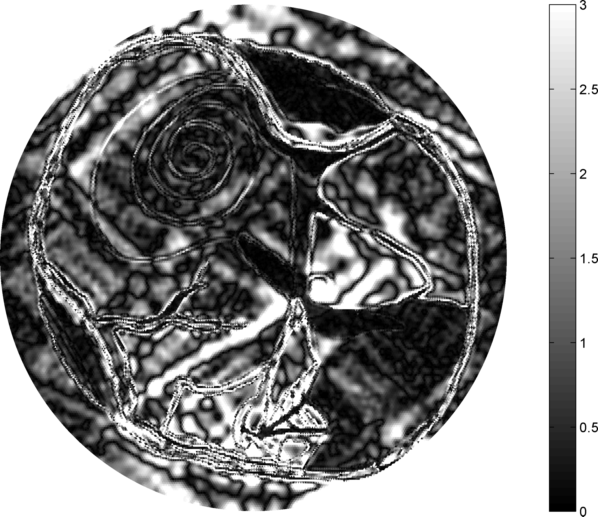}}  \hspace{0.01\textwidth}
  		\subfloat[$H^2$ adjoints, data $(\E_1^\delta,\E_3^\delta)$]{\includegraphics[width=0.32\textwidth]{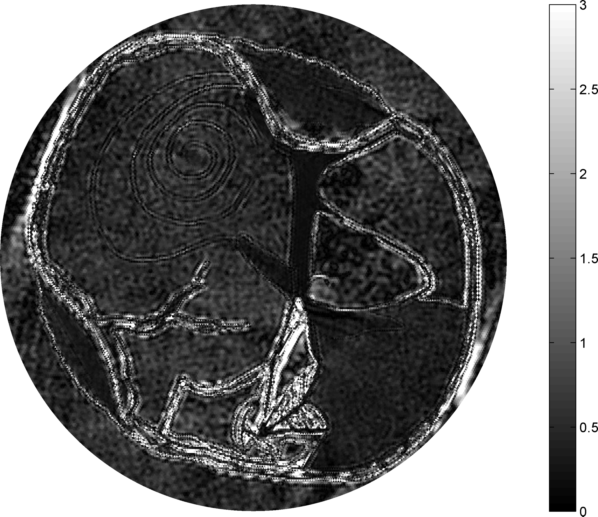}}  \hspace{0.01\textwidth}
  		\subfloat[$H^2$ adjoints, data $(\E_1^\delta,\E_2^\delta,\E_3^\delta)$]{\includegraphics[width=0.32\textwidth]{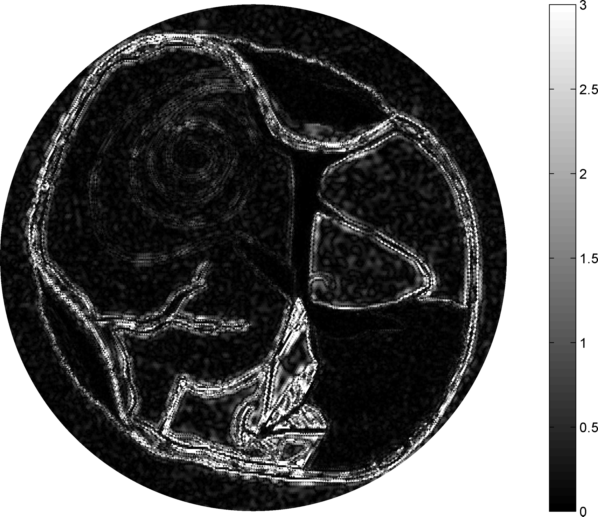}}

	\caption[Difference images]
	{
	  Difference images $|\sigma_{\text{est}}-\sigma|$ corresponding to the reconstructions in Figure \ref{fig:image_results}. The color axis was manually set to 0-3.
	}	
	\label{fig:image_difference}	
\end{figure}

The results of our numerical experiments seen in Figures \ref{fig:image_results} and \ref{fig:image_difference} show that from 2 measurements, even in the presence of significant noise, very good reconstructions are possible. The third measurement, which considerably increases the runtime, mostly serves to slightly reduce the noise. The $H^1$-approximation of the $H^2$ adjoint works well without loss of accuracy. The $L^2$-approximation, on the other hand, does not converge properly. From the difference images, we can see that in the main remaining error in the $H^1$ or $H^2$ reconstructions with 2 or more measurements is due to smoothing (which can be alleviated by running more iterations or lowering $\beta$).

While the given algorithms could in theory be directly translated to $\mathbb{R}^3$ (a higher degree of regularity would be necessary to get a Banach algebra and/or map into $L^2(\Omega)$), the drastically increased number of elements would make the computational effort for high-resolution reconstructions unreasonable. To reduce the effort, one could either switch to Landweber-like methods or use a combination of nested grids and iterative linear solvers as in \cite{CapFehGouKav09}.

\FloatBarrier

\section{Acknowledgements}
\label{sec:acknowledgements}

This work has been supported by the Austrian Science Fund (FWF) within the national research networks Photoacoustic Imaging in Biology and Medicine (project S10505) and Geometry+Simulation (project S11704) and by the IK I059-N funded by the University of Vienna. GB acknowledges partial support from the National Science Foundation Grant DMS-1108608.

%\bibliographystyle{plain}
%\bibliography{\BibPath strings,\BibPath articles,\BibPath books,\BibPath infmath,\BibPath infmath_books,\BibPath infmath_report,\BibPath infmath_talks,\BibPath infmath_theses,\BibPath inproceedings,\BibPath preprints,\BibPath proceedings,\BibPath theses,\BibPath unsubmitted,\BibPath websites}

\def\cprime{$'$} \providecommand{\noopsort}[1]{}

\end{document}